\documentclass[12pt]{amsart}
\usepackage{amsmath,amssymb,amsfonts,amscd,graphicx,xypic}

\usepackage{color}

\usepackage{latexsym}
\usepackage{amsmath}

\usepackage{xypic}
\usepackage[all]{xy}

\usepackage{hyperref}
\hypersetup{
    colorlinks = true,
   citecolor  = {blue},
}

\newtheorem{remark}[subsection]{Remark}
\numberwithin{equation}{section}

\newtheorem{cor}{Corollary}[section]
\newtheorem{prop}[cor]{Proposition}
\newtheorem{lem}[cor]{Lemma}
\newtheorem{theoremsection}[cor]{Theorem}
\newtheorem{conjecturesection}[cor]{Conjecture}
\theoremstyle{definition}

\theoremstyle{remark}

\newcommand{\Picture}[1]{
\begin{minipage}{.7in}
\includegraphics[scale=.2]{#1}
\end{minipage}
}

\newcommand{\PictureSmall}[1]{
\begin{minipage}{.6in}
\includegraphics[scale=.17]{#1}
\end{minipage}
}

\newcommand{\PictureSmalll}[1]{
\begin{minipage}{.5in}
\includegraphics[scale=.12]{#1}
\end{minipage}
}

\newcommand{\Ifig}{\ensuremath{[\!\Picture{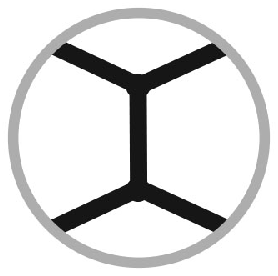}}\hspace{-2.5em}}
\newcommand{\Hfig}{\ensuremath{[\!\Picture{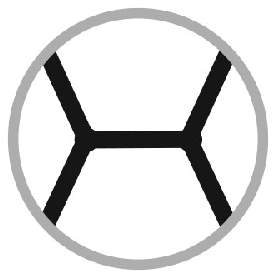}}\hspace{-2.5em}}
\newcommand{\ZeroRes}{\ensuremath{[\!\Picture{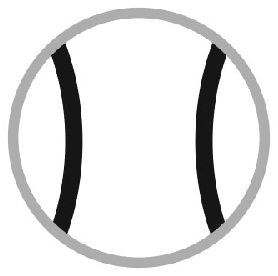}}\hspace{-2.5em}}
\newcommand{\OneRes}{\ensuremath{[\!\Picture{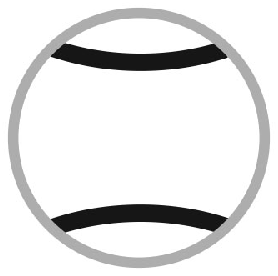}}\hspace{-2.5em}}
\newcommand{\Xfig}{\ensuremath{[\!\Picture{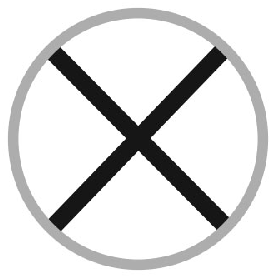}}\hspace{-2.5em}}
\newcommand{\Over}{\ensuremath{[\!\Picture{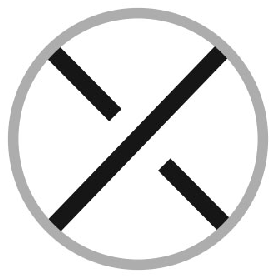}}\hspace{-2.5em}}
\newcommand{\Under}{\ensuremath{[\!\Picture{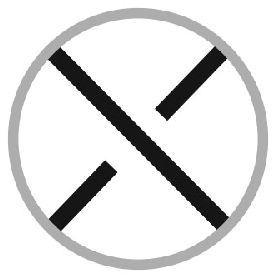}}\hspace{-2.5em}}

\newcommand{\ZeroResSmall}{\ensuremath{[\!\PictureSmall{0ResFig.eps}}\hspace{-2.5em}}
\newcommand{\OneResSmall}{\ensuremath{[\!\PictureSmall{1ResFig.eps}}\hspace{-2.5em}}
\newcommand{\XfigSmall}{\ensuremath{[\!\PictureSmall{Xfig.eps}}\hspace{-2.5em}}
\newcommand{\OverSmall}{\ensuremath{[\!\PictureSmall{OverCrossing.eps}}\hspace{-2.5em}}
\newcommand{\UnderSmall}{\ensuremath{[\!\PictureSmall{UnderCrossing.eps}}\hspace{-2.5em}}
\newcommand{\circlSmall}{\ensuremath{[\!\PictureSmall{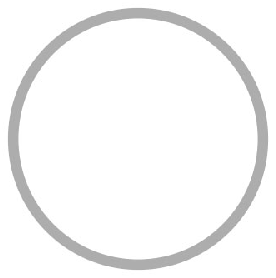}}\hspace{-2.5em}}

\newcommand{\circlblSmall}{\ensuremath{[\!\PictureSmalll{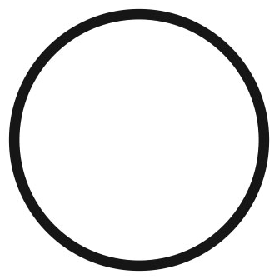}}\hspace{-2.5em}}

\begin{document}

\title[Flow and Yamada polynomials of cubic graphs]{Structure of the flow and Yamada polynomials of cubic graphs}

\author{Ian Agol and Vyacheslav Krushkal}

\address{Ian Agol{\hfil\break}
University of California, Berkeley, 970 Evans Hall \#3840, Berkeley, CA, 94720-3840}
\email{ianagol\char 64 berkeley.edu}

\address{Vyacheslav Krushkal{\hfil\break} Department of Mathematics, University of Virginia,
Charlottesville, VA 22904-4137 USA}
\email{krushkal\char 64 virginia.edu}

\begin{abstract} 
We establish a quadratic identity 
for the Yamada polynomial of ribbon  cubic graphs in ${\mathbb R}^3$, extending the Tutte golden identity for planar cubic graphs. An application is given to the structure of the flow polynomial of cubic graphs at zero.
The golden identity for the flow polynomial is conjectured to characterize planarity of cubic graphs, and we prove this conjecture for a certain infinite family of non-planar graphs.

Further, we establish exponential growth of the number of chromatic polynomials of planar triangulations, answering a question of D.~Treumann and E. Zaslow. 
The structure underlying these results is the chromatic algebra,  and more generally the ${\rm SO}(3)$ topological quantum field theory.
\end{abstract}

\maketitle

\section{Introduction}
Using the interplay between classical and quantum polynomials of graphs and ideas from topological quantum field theory (TQFT), we establish results on the structure of the Yamada and flow polynomials of cubic graphs. It has been known since the work of Birkhoff and Lewis in the 1940s \cite{BL} that the values $(3\pm \sqrt{5})/2$ of the parameter play a special role in the theory of the chromatic polynomial ${\chi}^{}_T$ of planar triangulations. In a series of papers \cite{T1, T2} in the 1960s, Tutte established further remarkable properties, including the  golden identity.  
Formulated dually in terms of the flow polynomial $F^{}_G$ of planar cubic graphs $G$, it reads 
\begin{equation} \label{flow golden identity}
F^{}_G({\phi}+2) \, =\,   {\phi}^{E}\, F^{}_G({\phi}+1)^2,
\end{equation}
where $E$ is the number of edges of $G$ and  ${\phi}$ denotes the golden ratio $(1+\sqrt{5})/2$. 
 The special role played by ${\phi}+1$, and more generally by the Beraha numbers \cite{Beraha}  $B_n=2+2{\rm cos}(2{\pi}/n)$, was conceptually explained in \cite{FK} where these results were placed in the context of ${\rm SO}(3)$ TQFT.
 
We show in Theorem \ref{Yamada golden} that Tutte's identity (\ref{flow golden identity}) admits an extension to the Yamada polynomial $R^{}_G$ of ribbon cubic graphs $G$ in ${\mathbb R}^3$: 
\begin{equation} \label{golden Yamada identity} 
R^{}_G(e^{{\pi} i/5})\, =\, (-1)^{V-E}\,  {\phi}^E\,  R^{}_G(e^{-2{\pi}i/5})^2,
\end{equation}
where $V, E$ denote the number of vertices and edges of $G$, respectively.
In fact, (\ref{golden Yamada identity}) is a common extension of (\ref{flow golden identity}) and of the identity for links, relating the $2$-colored Jones polynomial at $e^{{\pi} i/10}$ and the square of the Jones polynomial at $e^{-{\pi}i/5}$  \cite[Corollary 4.16]{MPS}, see section \ref{golden Yamada section}.

The Yamada polynomial \cite{Yamada} is a quantum invariant of ribbon graphs \cite{RT} in ${\mathbb R}^3$, corresponding to the adjoint representation of $U_q({\mathfrak{so}}_3)$.
Conceptually, as discussed in \cite[Section 5]{FK} and also \cite[Section 4.4]{MPS},  the reasons underlying the golden identity are the level-rank duality between the ${\rm SO}(3)_4$ and the ${\rm SO}(4)_3$ TQFTs, and the isomorphism ${\mathfrak{so}}(4)\cong {\mathfrak{so}}(3)\times {\mathfrak{so}}(3)$.

Concretely, the Yamada polynomial is defined by the contraction-deletion rule and the $SO(3)$ Kauffman skein relation,
 see section \ref{graph polynomials} for details.
For {\em planar} graphs $G$, the Yamada polynomial coincides with a renormalization of the flow polynomial:
$$F^{}_G(Q)\, =\, (-1)^{V-E}\, R^{}_G(q), \; {\rm where}\; Q=q+2+q^{-1}.$$
For non-planar graphs (and for knotted embeddings of planar graphs) the Yamada polynomial carries a lot of information about the embedding of a ribbon graph in $3$-space, and so in general the Yamada polynomial of a ribbon graph and the flow polynomial of the underlying abstract graph are quite different.

In contrast with (\ref{flow golden identity}), in \cite{AK} we formulated a conjecture that the Tutte golden identity for the flow polynomial characterizes planarity of cubic graphs:

\begin{conjecturesection} \label{golden conjecture0} \sl
For any cubic bridgeless graph $G$,
\begin{equation} \label{conjugate inequality0}
(-{\phi})^E\, F^{}_G((5-\sqrt{5})/2) \geq F^{}_G((3-\sqrt{5})/{2})^2,
  \end{equation}
Moreover, $G$ is planar if and only if (\ref{conjugate inequality0}) is an equality.
\end{conjecturesection}

The  inequality  at the Galois conjugate values,
\begin{equation} \label{conjugate conjecture}
F^{}_G({\phi}+2) \leq  {\phi}^E F^{}_G({\phi}+1)^2,
\end{equation}
conjecturally also holds for any cubic graph $G$, with an equality if and only if $G$ is planar.
In section \ref{conjecture section} we prove the conjecture for a family of near-planar graphs (which have a planar projection with a single crossing); Conjecture \ref{golden conjecture0} in general remains open.

It is interesting to note that as a consequence of (\ref{flow golden identity}), there is a relation between the values of the flow polynomial of planar cubic graphs at $0$ and $4$:
\begin{equation} \label{mod 5 F}
F^{}_G(0)\, \equiv \,  3^E\,  F^{}_G(4)^2 \;\,   ({\rm mod}\;  5),
\end{equation}
see lemma \ref{golden mod 5}. More generally, a congruence $({\rm mod} \, 5)$ between the values $R^{}_G(-1)$, $R^{}_G(1)$  implies an extension of (\ref{mod 5 F}) to {\em all} cubic graphs,
\begin{equation} \label{mod 5 Y}
F^{}_G(0)\, \equiv \,  3^E\,  R^{}_G(1)^2 \;\,   ({\rm mod}\;  5),
\end{equation}
where the value $R^{}_G(1)$ is also known as the Penrose number of $G$ (defined up to a sign), see section \ref{mod 5 section}.
Using these relations, 
we give an application to the structure  $({\rm mod} \, 5)$ of the flow polynomial at zero for cubic graphs. This value is known \cite{CS} to count Eulerian equivalence classes of totally cyclic orientations. 

\begin{theoremsection} \label{constant term theorem} \sl
Let $G$ be a cubic graph with $V$ vertices. Then the value of the flow polynomial of $G$ at zero satisfies: 
\begin{equation}  F^{}_G(0)\, \equiv 0,1,4 \; ({\rm mod} \, 5) \;\;  {\rm  if}\; \, V/2\; \, {\rm  is \;  even,}\end{equation}
\begin{equation}  F^{}_G(0) \, \equiv 0,2,3 \; ({\rm mod} \, 5)\; {\rm  if} \; \, V/2\; \,  {\rm  is \; odd}.
\end{equation}
Moreover, suppose $G$ is a snark, that is a bridgeless cubic graph with chromatic index $4$. Then $F^{}_G(0) \, \equiv 0\; ({\rm mod} \, 5)$, and therefore $F^{}_G(0)$ is divisible by $120$.
\end{theoremsection}

The proofs of the results in this paper are given in the context of the chromatic algebra \cite{FK}, and more generally  the ${\rm SO}(3)$ TQFT. Using these methods, we also  answer  a question of D. Treumann and E. Zaslow \cite{Treumann} about the asymptotics of the number of chromatic polynomials of planar triangulations, motivated by their work on Legendrian surfaces \cite{TZ}. 

\begin{theoremsection} \label{ChromaticPolynomials thm} \sl
The number of chromatic polynomials of planar triangulations with $n$ vertices grows exponentially in $n$.
\end{theoremsection}

An exponential upper bound is well-known, and it can be deduced from Tutte's enumeration of planar triangulations \cite{T63}. 
We show that the chromatic algebra contains a free semigroup, yielding an exponential lower bound, see section \ref{number of chromatic polynomials}. Conceptually, this may be viewed as an application of the Tits alternative for semigroups.

The key ingredients used throughout the paper -- the chromatic algebra and the flow category -- are recalled in section \ref{background}.
The identity (\ref{golden Yamada identity}) for the Yamada polynomial is established in section \ref{golden Yamada section}. Section  \ref{mod 5 section} discusses the structure of the flow polynomial at zero $ ({\rm mod} \, 5)$ and gives a proof of theorem \ref{constant term theorem}. The proof of conjecture \ref{golden conjecture} for a collection of near-planar graphs is given in section \ref{conjecture section}. Since the original Tutte golden identity (\ref{flow golden identity}) serves as the motivation for several results in this paper, for convenience of the reader we include its proof in section \ref{appendix}.

\section{Graph polynomials, algebras and categories} \label{background}
This section summarizes the relevant background material and notation used in the paper. 
\subsection{The chromatic and flow polynomials} \label{graph polynomials}
The flow polynomial $F^{}_G(Q)$ of a graph $G$ satisfies
the contraction-deletion rule: 
given an edge $e$ of $G$ which is not a bridge,  
\begin{equation}\label{contraction deletion rule}
F^{}_G(Q)\, =\, F^{}_{G/e}(Q)-F^{}_{G\smallsetminus e}(Q).
\end{equation}
If $G$ contains a bridge, $F^{}_G\equiv 0$.
The flow polynomial of a graph consisting of a single vertex and $n$ loops is defined to be $Q^n$, and the polynomial  is multiplicative with respect to taking the disjoint union. 

For planar graphs $G$, the flow polynomial is essentially the chromatic polynomial ${\chi}^{}_{G^*}$ of the dual graph ${G^*}$:  \begin{equation} \label{duality}
F^{}_G(Q)=Q^{-c} \, {\chi}^{}_{G^*}(Q),
\end{equation}
where $c$ is the number of connected components of $G$. The flow polynomial was defined by Tutte \cite{T3}, and in fact both the chromatic and flow polynomials are specializations of the $2$-variable Tutte polynomial $T^{}_G(x,y)$;  (\ref{duality}) is a special case of the duality relation satisfied by the Tutte polynomial: $T^{}_G(x,y)=T^{}_{G^*}(y,x)$.
Like the chromatic polynomial, the flow polynomial at positive integers admits a well-known combinatorial interpretation: for $n\in {\mathbb Z}_+$,  $F^{}_G(n)$ is the number of nowhere-zero flows with values in an abelian group of order $n$, cf. \cite{CS}.

\subsection{} \label{Yamada subsection} $\!$ {\bf The Yamada polynomial} \cite{Yamada} $R^{}_G(q)$ is an invariant of spatial ribbon graphs $G$, i.e. ribbon graphs embedded in ${\mathbb R}^3$. 
A {\em ribbon} graph is an abstract graph $G$ with an embedding into a
surface $S$, so that  the complement $S\smallsetminus G$ is a disjoint union of $2$-cells. Given such an embedding,
a neighborhood of $G$ in $S$ is a compact surface with boundary, which may be
thought of as a choice of a $2$-dimensional thickening of $G$. Such a thickening of vertices
and edges may be also encoded using cyclic ordering of half-edges incident to each
vertex.

The Yamada polynomial is defined using the contraction-deletion rule: 
\begin{equation}\label{contraction deletion rule1}
R^{}_G(Q)\, =\, R^{}_{G/e}(Q)+R^{}_{G\smallsetminus e}(Q), 
\end{equation}
and a version of the $SO(3)$ Kauffman skein relations (cf. \cite[Section 5]{FK2}) applied to a planar projection of $G$:
\begin{equation} \label{loop}
R^{}_{G\, \coprod\,  \circlblSmall}\; \, =\; (q+1+q^{-1})\, R^{}_G,
\end{equation}
\begin{equation} \label{over skein}
R^{}_{\OverSmall}\, =\, q\, R^{}_{\ZeroResSmall}\, +\, R^{}_{\XfigSmall}\, +\, q^{-1}\, R^{}_{\OneResSmall}
\end{equation}
\begin{equation} \label{under skein}
R^{}_{\UnderSmall}\, =\, q^{-1}\, R^{}_{\ZeroResSmall}\, +\, R^{}_{\XfigSmall}\, +\, q\, R^{}_{\OneResSmall}
\end{equation}
As usual in skein-theoretic definitions, the graphs in each of these equations differ as shown, and are identical outside of the disk. If $G$ has a bridge, $R^{}_G$ is set to be zero.
The Yamada polynomial is multiplicative under taking disjoint union, and $ R^{}_{G\, \vee\, H}\, =\, - R^{}_{G}R^{}_{H}$.
 $R^{}_G$ is an invariant of the isotopy class of an embedding of the ribbon graph $G$ into ${\mathbb R}^3$. (In terminology of \cite[Theorem 2]{Yamada},
$R^{}_G$ is a {\em regular deformation} invariant  of a planar diagram of $G$; in \cite{KV} this equivalence relation is called {\em rigid vertex isotopy}.)

Using just equations (\ref{loop}) -- (\ref{under skein}), one gets a (renormalized version of) the SO$(3)$ Kauffman polynomial of $4$-regular ribbon graphs in ${\mathbb R}^3$ (cf. \cite{KV}, \cite{FK2}). Therefore the Yamada polynomial may be thought of as the SO$(3)$ Kauffman polynomial, extended to spatial ribbon graphs of arbitrary vertex degree using the contraction-deletion rule (\ref{contraction deletion rule1}).

If $G$ is {\em planar}, then there are no crossings to resolve, and the Yamada polynomial is determined by the contraction-deletion rule and its loop value (\ref{loop}); in this case $F^{}_G(Q)\, =\, (-1)^{V-E}\, R^{}_G(q)$,  where $Q=q+2+q^{-1}.$ In general, the Yamada polynomial carries a lot of information about the embedding of $G$ into ${\mathbb R}^3$. For example, tying a knot into an edge of $G$ results (up to a normalization) in multiplication of $R^{}_G$ by the SO$(3)$ invariant of the knot. Therefore, in general the Yamada polynomial of a ribbon graph is quite different from the flow polynomial of the underlying abstract graph.

To describe the TQFT context for the results of this paper, 
next we give a brief summary of the relevant material on the Temperley-Lieb algebra, the chromatic algebra and their structure at roots of unity.

\subsection{The Temperley-Lieb algebra} \label{TL subsection}
\begin{figure}[ht]
\centering
\includegraphics[height=1.7cm]{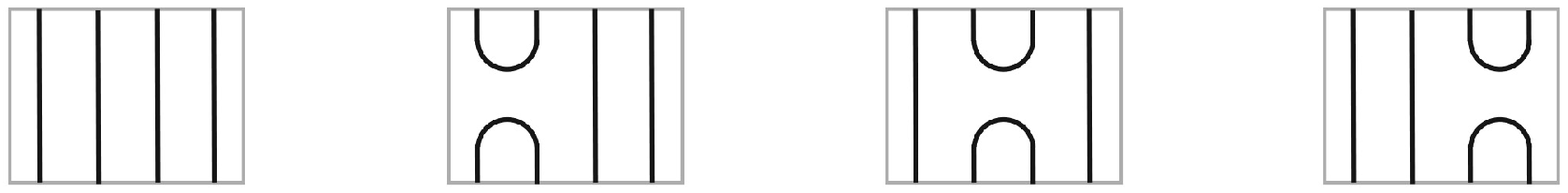}
{\small     \put(-350,20){$1 =$}
    \put(-271,20){$e_1 =\frac{1}{d}$}
     \put(-180,20){$e_2 = \frac{1}{d}$}
    \put(-89,20){$e_3 = \frac{1}{d}$}
    }  
    \caption{Generators of ${\rm TL}_4$}
\label{TL figure}
\end{figure}

The Temperley-Lieb algebra, ${\rm TL}_n$, is an algebra over
${\mathbb C}[d]$ consisting of linear combinations of
$1$-dimensional submanifolds, considered up to isotopy rel boundary, in a rectangle. 
Each submanifold meets both the top and the bottom of the rectangle in $n$ points. 
Deleting a simple closed curve has the effect of multiplying the element by $d$. 
Often $d$ will be specialized to a complex number, and in this case the algebra will be denoted
${\rm TL}^d_n$.
The multiplication is given by vertical stacking of rectangles. The standard generators
of ${\rm TL}_4$ are shown in figure \ref{TL figure}.

The trace ${\rm tr}_d{\colon\thinspace} {\rm TL}^d_n\longrightarrow {\mathbb C}$ is defined on
rectangular pictures by connecting the top
and bottom endpoints by disjoint arcs in the complement of the rectangle in the
plane,
and then evaluating $d^{\# {\rm circles}}$. The Hermitian product 
is defined by $\langle a,b\rangle={\rm tr}(a\, \overline b)$, where the
involution $^-$ reflects pictures in a horizontal line and
replaces coefficients with their complex conjugates.

For special values of the parameter $d$, $d\, = \, 2\, \cos \left(\frac{\pi k}{n+1}\right)$, where $k, n+1$ are coprime, ${\rm TL}^d_n$ contains a non-trivial ideal, the {\em trace radical} consisting  of
the elements $a$ such that $tr(ab)=0$ for all $b\in {\rm TL}^d_n$.
This ideal is generated by the {\em Jones-Wenzl
projector}  $P^{(n)}$ \cite{We, Jo}.
At primitive roots of unity ($k=1$) the Hermitian product descends to a 
positive definite inner product on the quotient of ${\rm TL}^d_n$ by the trace radical.

\subsection{The chromatic algebra ${\mathbf{{\mathcal C}_n}}$} \label{chromatic algebra subsection}
The Temperley-Lieb algebra, discussed above, underlies the construction of SU$(2)$ TQFT.  Next we briefly summarize the definition and properties of the chromatic algebra, ${\mathcal C}_n$, corresponding to SO$(3)$ TQFT; we refer the reader to \cite[Section 2]{FK} for more details. (A similar notion, in a different context, was considered in \cite[Section 2]{MW}.)

\begin{figure}[ht]
\includegraphics[height=1.7cm]{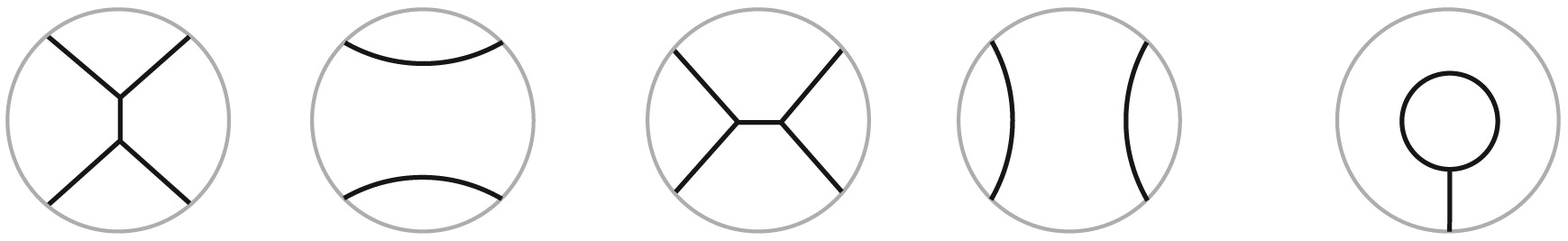}
    \put(-269,20){$+$}
    \put(-206,20){$=$}
    \put(-139,20){$+$}
 {\small    \put(2,20){$=0.$}}
    \put(-79,8){$,$}
\caption{Relations defining the chromatic algebra.}
\label{chromatic relations}
\end{figure}

${\mathcal C}_n$ is an algebra over ${\mathbb C}[Q]$, 
whose elements are formal linear combinations of 
 {\em planar} cubic graphs (considered up to isotopy rel boundary) in a
rectangle, modulo the local relations shown in figure \ref{chromatic relations}.
The first relation is the analogue  of the contraction-deletion relation for cubic graphs. 

\begin{figure}[ht]
\includegraphics[height=1.7cm]{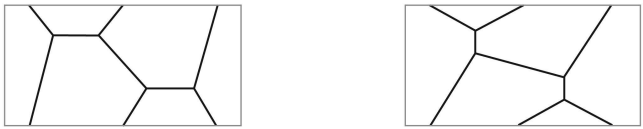}
\caption{Examples of graphs in ${\mathcal C}_3$.}
\label{Chromatic Generator}
\end{figure}

The intersection of a graph with the boundary
of the rectangle consists of $2n$ points: $n$ points both at the top and
the bottom, figure \ref{Chromatic Generator}. 
It is convenient to allow $2$-valent vertices as well, 
and the value of a
simple closed curve is set to be $Q-1$. When $Q$ is specified to a complex number, 
the algebra is denoted ${\mathcal C}^Q_n$.

The trace $tr{\colon\thinspace} {\mathcal C}_n^Q\longrightarrow{\mathbb C}$ is defined by connecting the
endpoints of $G$ by disjoint arcs in the complement of the rectangle in the plane
and evaluating the flow polynomial of the resulting graph at $Q$. (Or equivalently, the trace equals $Q^{-1}$ times the chromatic 
polynomial of the dual graph.)
The trace is well-defined since the local relations in figure \ref{chromatic relations} are precisely the relations defining the flow polynomial of a planar cubic graph.
The multiplication and the Hermitian product on ${\mathcal C}^Q_n$ are defined analogously to those in the Temperley-Lieb algebra.

There are two variations of the definition of the chromatic algebra, which are going to be useful.
First, rather than using just cubic graphs modulo relations in figure \ref{chromatic relations}, ${\mathcal C}_n$ may also be defined using graphs with arbitrary vertex degrees, modulo the contraction-deletion relation, see \cite[Section 5]{FK}.

Second, instead of using planar graphs, one may consider ribbon graphs in the cylinder $D^2\times [0,1]$, with $n$ endpoints both at $D^2\times 0$ and at $D^2\times 1$, modulo the defining relations of the Yamada polynomial, (\ref{contraction deletion rule1}) - (\ref{under skein}). This is closely related to the definition of the SO$(3)$ BMW algebra (cf. \cite[Section 5]{FK2}).
Using the relations $${\Over}\, =\, q\, {\ZeroRes}\, +\, {\Xfig}\, +\, q^{-1}\,{\OneRes}, \; \, {\UnderSmall}\, =\, q^{-1}\, {\ZeroRes}\, +\, {\Xfig}\, +\, q\, {\OneRes},$$
all crossings of a ribbon graph  in the cylinder may be resolved to give an element of ${\mathcal C}_n$.

\subsection{The map ${\mathbf{{\mathcal C}_n\longrightarrow {\rm TL}_{2n}}}$} \label{map}

\begin{figure}[ht]
\includegraphics[height=1.8cm]{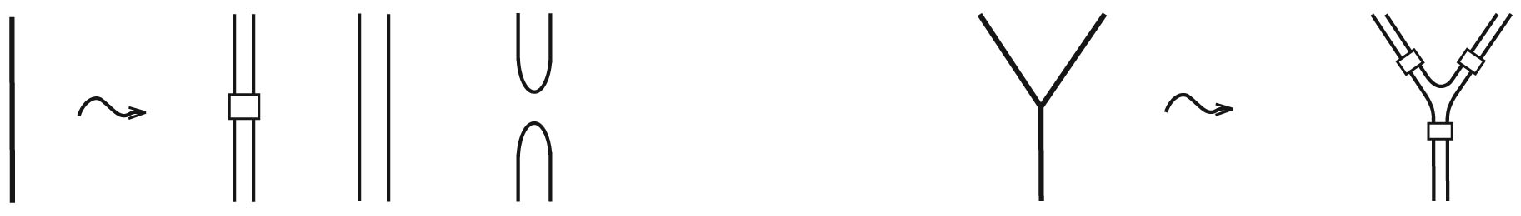}
    \put(-265,21){$=$}
    \put(-239,21){$-\frac{1}{d}$}
    \put(-54,21){$d^{1/2}\;\cdot$}
\caption{The homomorphism ${\Phi}{\colon\thinspace} {\mathcal C}^Q_n\longrightarrow {\rm TL}^d_{2n}$, $Q=d^2$.}
\label{map fig}
\end{figure}

Consider the homomorphism ${\Phi}{\colon\thinspace} {\mathcal C}^Q_n\longrightarrow
{\rm TL}^d_{2n}$, where $Q=d^2$,  replacing 
each  edge of a graph  with the second Jones-Wenzl projector,  
and resolving  
each vertex as shown in figure \ref{map fig}. Moreover, for a trivalent graph $G$ there is an overall factor
$d^{V/2}$, where $V$ is the number of vertices of $G$.

${\Phi}$ induces a well-defined homomorphism of algebras ${\mathcal C}^Q_n\longrightarrow {\rm TL}^d_{2n}$,
where $Q=d^2$, and moreover it is trace-preserving: the diagram
\begin{equation} \label{chromaticTL}
\xymatrix{ {\mathcal C}_n^{{d^2}}  \ar[d]^{tr} \ar[r]^{\Phi} & {\rm TL}_{2n}^{d}  \ar[d]^{tr_{d}}\\
{\mathbb C}\ar[r]^=  &  {\mathbb C} }\end{equation}
commutes \cite[Lemmas 2.4, 2.5]{FK}.
It follows that the pullback under ${\Phi}$
of the trace radical in ${\rm TL}^d$ is in the trace radical of ${\mathcal
C}^{d^2}$. As mentioned in section \ref{TL subsection}, the trace radical in ${\rm TL}^d$ is non-trivial precisely for the values $d\, = \, 2\, \cos \left(\frac{\pi k}{n+1}\right)$. The elements of the trace radical  in the chromatic algebra ${\mathcal
C}^{d^2}$
are local relations (which hold in addition to the contraction-deletion rule) on graphs which preserve the flow polynomial, 
or equivalently the chromatic polynomial of the dual graph.
When $d={\phi}$, the trace radical of the Temperley-Lieb algebra is
generated by the Jones-Wenzl projector $P^{(4)}$.   In particular, the linear Tutte relation
in ${\mathcal C}^{Q}_4$, where $Q={\phi}+1={\phi}^2$,
\begin{equation} \label{linear Tutte0}
\Ifig\, =\, {\phi}^{-1}\, \ZeroRes - {\phi}^{-2}\,  \OneRes,
\end{equation}

may be seen as a consequence of the structure of ${\rm TL}^{\phi}$ since it is mapped by $\Phi$ to
$P^{(4)}$ \cite[Section 2, p. 721]{FK}. Similarly, the linear relation
\begin{equation} \label{linear Tutte1}
\Ifig\, =\, -{\phi}\, \ZeroRes - {\phi}^{2}\,  \OneRes,
\end{equation}

holds in ${\mathcal C}^{Q}_4$, at the Galois conjugate value $Q=(3-\sqrt{5})/2$.

\subsection{The flow category} \label{flow category subsection}
To  study the flow polynomial of non-planar graphs, one may use abstract (not necessarily planar) graphs to define  an algebra 
along the lines of the Temperley-Lieb and chromatic algebras \cite{MW, Walker}. 
Unlike the chromatic algebra case, here one does not have a map to the Temperley-Lieb algebra, and the values $Q=d^2, \, d\, = \, 2\, \cos \left(\frac {\pi}{n+1}\right)$ do  not play a special role since these are artifacts of planarity.
We will not use the algebra structure; the notion of a category is more suitable for our applications. The construction of the {\em flow category} is summarized below, following  \cite[Section 3.1]{AK}. 

The objects of the flow category are finite ordered sets $\overline{n}=\{ 1,\ldots, n\}$. 
Consider ${\mathcal G}^{}_{m,n}=\{$finite graphs with $m+n$  marked univalent vertices$\}$, where the marked vertices are divided into two ordered subsets of $m$, respectively $n$ vertices. The edges incident to the marked vertices
are called  boundary edges and the rest are {\em internal} edges. 

\begin{figure}[h]
\begin{center}
\includegraphics[height=1.9cm]{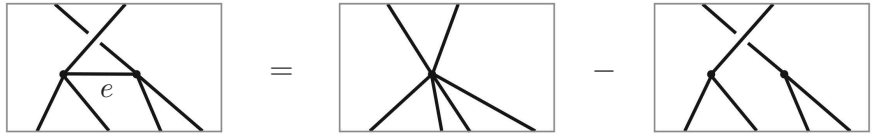}
\put(-230,23){$=$}
\put(-108,23){$-$}
{\small \put(-294,16){$e$}}
\caption{The contraction-deletion rule in ${\mathcal F}_{4,2}$.}
\label{Flow category}
\end{center}
\end{figure}

The space of morphisms ${\mathcal F}^Q_{m,n}$ in the flow category between $\overline{m}$, $\overline{n}$  consists of  formal ${\mathbb C}$-linear
combinations of graphs ${\mathcal G}^{}_{m,n}$, modulo the contraction-deletion relation which applies to internal edges,  figure \ref{Flow category}.
Graphs whose equivalence classes are elements of ${\mathcal F}^Q_{m,n}$ may be represented geometrically as in figure \ref{Flow category}.
(It is important to note that unlike in sections \ref{Yamada subsection}, \ref{chromatic algebra subsection}, over/under-crossings do not carry any information here since the figure represents the abstract graph structure and not a specific planar projection.)
In addition, the loop value is set to be $Q-1$, and graphs with a univalent vertex (other than the specified marked vertices) 
are set to be zero. 

A graph without marked vertices, and therefore no boundary edges, considered in ${\mathcal F}^Q_{0,0}\cong {\mathbb C}$, evaluates to its flow polynomial at $Q$.
The pairing ${\mathcal F}^Q_{k,m}\times {\mathcal F}^Q_{m,n}\longrightarrow {\mathcal F}^Q_{k,n}$ is obtained by gluing along $m$  boundary edges.
For example, this pairing applied to two graphs $A\in {\mathcal G}^{}_{0,m}$, $B\in {\mathcal G}^{}_{m,0}$ gives $\langle A, B\rangle=$ the value of the flow polynomial 
$F^{}_{A\cup B}(Q)$.

\begin{figure}[h]
\begin{center}
\includegraphics[height=1.75cm]{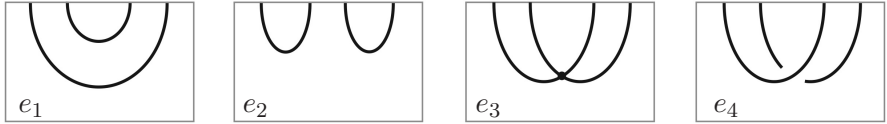}
{\small    \put(-330,5){$e_1$}
    \put(-245,5){$e_2$}
    \put(-157,5){$e_3$}
    \put(-67,5){$e_4$}
    }
\caption{A basis of ${\mathcal F}^Q_{0,4}$.}
\label{basis figure}
\end{center}
\end{figure}

Given any graph representing an element of ${\mathcal F}^Q_{m,n}$,  the contraction-deletion rule may be used to eliminate all internal edges. For example, the four graphs  in figure \ref{basis figure} form a basis of  ${\mathcal F}^Q_{0,4}$. Three of these graphs, viewed relative to a fixed embedding of the marked vertices in the boundary, are planar, and one, denoted $e_4$, is non-planar.

It is convenient to introduce the notation ${\mathcal C}^Q_{m,n}$ for the planar analogue of ${\mathcal F}^Q_{m,n}$. That is, ${\mathcal C}^Q_{m,n}$  consists of  formal ${\mathbb C}$-linear
combinations of {\em planar} graphs with two ordered subsets of $m$, respectively $n$ marked vertices, modulo the contraction-deletion relation which applies to internal edges.

\subsection{Conventions and notation.} \label{notation section}

Unless stated otherwise,  in the following sections the flow polynomial of {\em planar} graphs and the Yamada polynomial of spatial ribbon graphs will be considered in the context of the chromatic algebra. On the other hand, the flow polynomial of abstract (non-planar) graphs does not fit in this context, and it will be studied in the setting of the flow category.

It is convenient to introduce a short-hand notation for the evaluation of graph polynomials. For example, when working with the Yamada polynomial of a graph $G$ with a specified crossing $\Over$, the notation $\Over_{\!\! x}$ will stand for the evaluation of $R^{}_G(x)$. Similarly, given two graphs $G_1, G_2$ in a disk with the same number of marked points on the boundary, the notation $\langle G_1, G_2\rangle_Q$ will stand for the relevant pairing. For example, in the context of the flow category, it will mean the value of the flow polynomial $F^{}_{G_1\cup G_2}$ at $Q$, where the union of $G_1, G_2$ is taken along the marked points on the boundary.

\section{The golden identity for the Yamada polynomial} \label{golden Yamada section}

The purpose of this section is to prove the extension (\ref{golden Yamada identity}) of the Tutte golden identity (\ref{flow golden identity}) to the Yamada polynomial of cubic ribbon graphs in ${\mathbb R}^3$. It is convenient to allow vertices of degree $2$ in the statement of the theorem:

\begin{theoremsection} \label{Yamada golden} \sl
Let $G$ be a ribbon graph in ${\mathbb R}^3$, with vertices of degrees $2$ and $3$. Then 
\begin{equation} \label{golden Yamada} 
R^{}_G(e^{{\pi} i/5})\, =\, (-1)^{V-E}\,  {\phi}^{E'}\,  R^{}_G(e^{-2{\pi}i/5})^2,
\end{equation}

where $E'=V_3-{\chi}(G)$,   $V_3$ is the number of trivalent vertices of $G$, ${\chi}(G)$ is the Euler characteristic, and $\phi=(1+\sqrt{5})/2$.
\end{theoremsection}

If there are no vertices of degree $2$, then $E'=E$, and (\ref{golden Yamada}) is the same identity as (\ref{golden Yamada identity}). Moreover, since both ${\chi}(G)$ and $V_3(G)$ are topological invariants of $G$, introducing new $2$-valent vertices (i.e. subdividing edges of $G$) does not affect (\ref{golden Yamada}).

In the special case where $G$ has no trivalent vertices, (\ref{golden Yamada}) gives an identity for  the SO$(3)$ Kauffman polynomial of framed {\em links} in ${\mathbb R}^3$, where the factor $(-1)^{V-E}\,  {\phi}^{E'}$ equals $1$. Then the left-hand side of  (\ref{golden Yamada}) may be interpreted as the $2$-colored Jones polynomial of the link, and the right-hand side equals the square of the Jones polynomial. Therefore for links, (\ref{golden Yamada}) matches the identity stated in \cite[Corollary 4.16]{MPS}.

{\em Proof of theorem \ref{Yamada golden}.}
Recall from section \ref{Yamada subsection}
that if $G$ is {\em planar}, one has $F^{}_G(Q)\, =\, (-1)^{V-E}\, R^{}_G(q)$,  where $Q=q+2+q^{-1}.$ In particular, for planar graphs $G$, 
$$R^{}_G(e^{2{\pi} i/10})\, =\, F^{}_G({\phi}+2), \; \, R^{}_G(e^{-2{\pi} i/5})\, =\, F^{}_G({\phi}+1).$$      

Therefore in this case, (\ref{golden Yamada}) is reduced to the Tutte golden identity (\ref{flow golden identity}).\footnote{Its proof is given in section \ref{appendix}.}
In the general case of a ribbon graph $G$ in the statement of the theorem, the proof is by induction on the number $c$ of crossings in a planar diagram of $G$. The base case corresponds to  planar graphs, discussed above.

As in section \ref{notation section}, 
the notation $\circlSmall_{\; q}$ in the proof below will denote the evaluation of the Yamada polynomial at $q$. 

By induction assume that graphs with fewer than $c$ crossings satisfy (\ref{golden Yamada}). Consider a graph with $c$ crossings. For brevity of notation denote $x:=e^{2\pi i/10}, y:=e^{-2\pi i/5}$. 
Combining the skein relation (\ref{over skein}) with the contraction-deletion rule (\ref{contraction deletion rule1}), one has 
\begin{equation} \label{calc1}
{\Over}_{\!\! x} = e^{2{\pi} i/10}\, {\ZeroRes}_{\!\! x} + {\Ifig}_{\!\! x} +(e^{-2\pi i/10}-1) \, {\OneRes}_{\!\! x}.
\end{equation}

Since the three graphs on the right have $c-1$ crossings, (\ref{golden Yamada}) holds for them by the inductive assumption, thus
\begin{equation} \label{calc2}
{\Over}_{\!\! x} = (-1)^{V-E} {\phi}^{E'}\left ( e^{2{\pi} i/10}\, {\ZeroRes}_{\!\! y}^{\! 2} - {\phi}^3 {\Ifig}_{\!\! y}^{\! 2} +(e^{-2\pi i/10}-1) \, {\OneRes}_{\!\! y}^{\! 2} \right ).
\end{equation}

Because of the normalization sign $(-1)^{V-E}$ in the formula relating the flow and Yamada polynomials of planar graphs, the linear relation (\ref{linear Tutte0}) at $q=e^{-2 \pi i/5}$ (corresponding to $Q={\phi}+1$) for the Yamada polynomial reads
\begin{equation} \label{linear Tutte00}
{\Ifig}_{\!\! y}\, =\, - {\phi}^{-1}\, {\ZeroRes}_{\!\! y} + {\phi}^{-2}\,  {\OneRes}_{\!\! y},
\end{equation}

Applying this relation to ${\Ifig}_{\!\! y}$ in (\ref{calc2}) gives
\begin{equation} \label{calc4}
{\Over}_{\!\! x} = (-1)^{V-E}{\phi}^{E'}\left ( (e^{2{\pi} i/10}-{\phi}) \, {\ZeroRes}_{\!\! y}^{\! 2} +2\, {\ZeroRes}_{\!\! y} \, {\OneRes}_{\!\! y} +(e^{-2\pi i/10}-{\phi}) \, {\OneRes}_{\!\! y}^{\! 2} \right ).
\end{equation}

To prove the inductive step, one needs to show that (\ref{calc4}) equals $(-1)^{V-E}  {\phi}^{E'} {\Over}_{\!\! y}^{\! 2}$. Using the skein relation  (\ref{over skein}) and the contraction-deletion rule, 
\begin{equation} \label{calc5}
(-1)^{V-E}  {\phi}^{E'} {\Over}_{\!\! y}^{\! 2} = (-1)^{V-E} {\phi}^{E'}\left ( e^{-2{\pi} i/5}\, {\ZeroRes}_{\!\! y} + {\Ifig}_{\!\! y} +(e^{2\pi i/5}-1) \, {\OneRes}_{\!\! y} \right )^2.
\end{equation}
To complete the proof, one replaces ${\Ifig}_{\!\! y}$ with (\ref{linear Tutte00}), and checks that the resulting expression matches  (\ref{calc4}).
\qed

\section{Structure of the flow polynomial (mod $5$)} \label{mod 5 section}
In this section we establish the $({\rm mod}\;  5)$ version of the golden identity, and prove a theorem stated in the introduction:

\medskip

{\bf Theorem 1.2.} {\sl
Let $G$ be a cubic graph with $V$ vertices. Then the value of the flow polynomial of $G$ at zero satisfies: 
\begin{equation} \label{theorem eq 1} F^{}_G(0)\, \equiv 0,1,4 \; ({\rm mod} \, 5) \;\;  {\rm  if}\; \, V/2\; \, {\rm  is \;  even,}\end{equation}
\begin{equation} \label{theorem eq 2} F^{}_G(0) \, \equiv 0,2,3 \; ({\rm mod} \, 5)\; {\rm  if} \; \, V/2\; \,  {\rm  is \; odd}.
\end{equation}
Moreover, suppose $G$ is a snark, that is a bridgeless cubic graph with chromatic index $4$. Then $F^{}_G(0) \, \equiv 0\; ({\rm mod} \, 5)$, and therefore $F^{}_G(0)$ is divisible by $120$.}

{\bf  Remark.} \cite[Theorem 1.2]{CS} interpreted $|F^{}_G(0)|$ as the number of Eulerian equivalence classes of totally cyclic orientations.

{\em Proof of theorem \ref{constant term theorem}.} 
We start  by stating the $({\rm mod}\;  5)$ version of the golden identity:

\begin{lem} \label{golden mod 5} \sl
For planar cubic graphs $G$, 
\begin{equation} \label{mod 5 flow eq} 
F^{}_G(0)\, \equiv \,  3^E\,  F^{}_G(4)^2 \;\,   ({\rm mod}\;  5),
\end{equation}
where $E$ is the number of edges of $G$. More generally, let $G$ be a ribbon cubic graph in ${\mathbb R}^3$. Then 
\begin{equation} \label{mod 5 Yamada eq} 
R^{}_G(-1)\, \equiv \,  (-1)^{V-E}\, 3^E\,  R^{}_G(1)^2 \;\,   ({\rm mod}\;  5).
\end{equation}
\end{lem} 

{\em Proof of} (\ref{mod 5 flow eq}). The flow polynomial $F^{}_G(Q)$ for $Q=(3-\sqrt{5})/2, (5-\sqrt{5})/2$ takes values in the ring $R={\mathbb Z}[\frac{1+\sqrt{5}}{2}]$.
Consider the ideal $I$ generated by $\sqrt{5}$. Since ${\phi}\sqrt{5}=(5+\sqrt{5})/2 \in I$, it follows that
$  1+{\phi} \equiv 4\;  ({\rm mod} \; \sqrt{5})$, and  $ {\phi} \equiv 3\;  ({\rm mod} \; \sqrt{5})$. Plugging this into the golden identity (\ref{flow golden identity}), one gets
$F^{}_G(0)  \equiv 3^E F^{}_G(4)^2 \; ({\rm mod} \; \sqrt{5})$. 
Since $ F^{}_G(0)$ and $ 3^E F^{}_G(4)^2$ are integers, the equivalence holds  $({\rm mod} \; 5)$.
\qed

{\em Proof of} (\ref{mod 5 Yamada eq}). The values of the Yamada polynomial $R^{}_G(e^{ {\pi} i/5}),  R^{}_G(e^{-2{\pi}i/5})$ are elements of ${\mathbb Z}[{\zeta}_5]$, where ${\zeta}_5=e^{{\pi}i/5}$. Consider the ideal generated by $e^{{\pi} i/5}+1$. Modulo this ideal, $e^{{\pi} i/5}$ and $e^{-{\pi} i/5}$ are equivalent to  $-1$, $e^{2{\pi} i/5}$ and $e^{-2{\pi} i/5}$ are equivalent to $1$, and ${\phi}=e^{{\pi} i/5}+e^{-{\pi} i/5}\equiv 3$.
\qed

\smallskip

{\bf Remark.}
Another proof of lemma \ref{golden mod 5} may be given  by following $({\rm mod}\;  5)$ the steps of the proofs of the golden identities (\ref{flow golden identity}), (\ref{golden Yamada}). For example, the (mod $5$) version at $Q=4$ of the Tutte linear relation (\ref{linear Tutte0}) for the flow polynomial of planar graphs, ${\Ifig}_{\!\! {\phi}+1} = {\phi}^{-1}{\ZeroRes}_{\!\!{\phi}+1} -{\phi}^{-2} {\OneRes}_{\!\! {\phi}+1}$, reads ${\Ifig}_{\!\! 4}\, \equiv \, 2\; {\ZeroRes}_{\!\! 4} + {\OneRes}_{\!\! 4}$ $({\rm mod} \; 5)$.

\smallskip

Returning to the proof of theorem \ref{constant term theorem},
recall that for planar graphs $G$, $F^{}_G(Q)\, =\, (-1)^{V-E}\, R^{}_G(q)$,  where $Q=q+2+q^{-1}.$ In particular, for planar $G$, $R^{}_G(-1)\, =\,  (-1)^{V-E} F^{}_G(0)$.
The following identity holds for the flow polynomial of all (planar and non-planar) graphs:
\begin{equation} \label{flow at 0 eq}
{\Over}_{\!\! 0}\; = \; - {\ZeroRes}_{\!\! 0} - {\Xfig}_{\!\! 0} - {\OneRes}_{\!\! 0}.
\end{equation}
This identity is checked by pairing both sides of (\ref{flow at 0 eq}) with the four basis elements 
of the flow category ${\mathcal F}^Q_{0,4}$ in figure \ref{basis figure}. This relation coincides with the skein relation  (\ref{over skein}) for the Yamada polynomial,  normalized by the factor $(-1)^{V-E}$,  at $q=-1$. Thus $R^{}_G(-1)\, =\,  (-1)^{V-E} F^{}_G(0)$ for {\em all} graphs $G$.
Now it follows from (\ref{mod 5 Yamada eq}) that for any abstract (planar or non-planar) cubic graph $G$, 
\begin{equation} \label{Penrose number squared}
F^{}_G(0)\, \equiv \,   3^E\,  R^{}_G(1)^2 \;\,   ({\rm mod}\;  5).
\end{equation}

{\bf Remark.} Up to a sign, the value $R^{}_G(1)$ of the Yamada polynomial
of ribbon cubic graphs equals the {\em Penrose number} of $G$, introduced and studied in \cite{Penrose}, \cite[Section 2.3]{Jaeger}.  Indeed, the skein relation \cite[Proposition 2]{Jaeger}, satisfied by the Penrose number, coincides with the version of the skein relation (\ref{over skein}) for $(-1)^{V-E}R^{}_G(1)$:
\begin{equation}\label{skein at 1}
 {\Ifig}_{\!\! 1}   = {\ZeroRes}_{\!\! 1} - {\Over}_{\!\! 1}.
 \end{equation}
 It is an invariant of ribbon graphs where the ribbon structure affects only its sign, so $|R_G(1)|$ is an invariant of abstract cubic graphs. 

The congruence (\ref {Penrose number squared}), the relation $3V=2E$ for cubic graphs, and the fact that perfect squares are congruent to $0, 1$, or $4$ $({\rm mod}\; 5)$  conclude the proof of (\ref{theorem eq 1}), (\ref{theorem eq 2}).

To prove the last statement of the theorem, consider a planar diagram of a ribbon cubic graph $G$. 
In other words, the graph is immersed in the plane, with some
edge crossings. The value $R^{}_G(1)$ is independent of which strand is over/under in each crossing.
Using the recursion relation (\ref{skein at 1}) at $q=1$, one checks that $(-1)^{V-E} R^{}_G(1)$ equals the signed count of $3$-edge colorings of $G$,
where the sign gets a $(-1)$ factor for each pair of edges that
cross with a different color. This is invariant
of regular homotopy of the planar diagram, and it satisfies the skein relation: each $3$-coloring of the graph on the left in (\ref{skein at 1}) corresponds (with the same sign) to a $3$-coloring of precisely one term on the right.

This implies that the Yamada polynomial of snarks at $q=1$ is zero. (More generally, there is a well-known correspondence between $3$-edge-colorings and $4$-flows, cf. Proposition 6.4.5 (ii) in \cite{Diestel}, so for any cubic graph $|R^{}_G(1)|\leq F^{}_G(4)$, with equality
for planar graphs.)
Now the congruence $F^{}_G(0) \, \equiv 0\; ({\rm mod} \, 5)$ follows from (\ref{Penrose number squared}). 
Finally, since the flow polynomial $F^{}_G(Q)$ of a snark is divisible by $\prod_{k=1}^4 (Q-k)$, $F^{}_G(0)$ is also divisible by $4!$. \qed

\section{Golden inequality for non-planar cubic graphs} \label{conjecture section}

Unlike the golden identity for the Yamada polynomial, Theorem \ref{Yamada golden}, which holds for any spatial cubic graph $G$, in \cite{AK} we stated a conjecture that the golden identity for the flow polynomial characterizes planarity:

\begin{conjecturesection} \label{golden conjecture} \sl
For any cubic bridgeless graph $G$,
\begin{equation} \label{conjugate inequality}
(-{\phi})^E\, F^{}_G((5-\sqrt{5})/2) \geq F^{}_G((3-\sqrt{5})/{2})^2,
  \end{equation}
Moreover, $G$ is planar if and only if (\ref{conjugate inequality}) is an equality.
\end{conjecturesection}

\begin{figure}[ht]
\includegraphics[height=2.5cm]{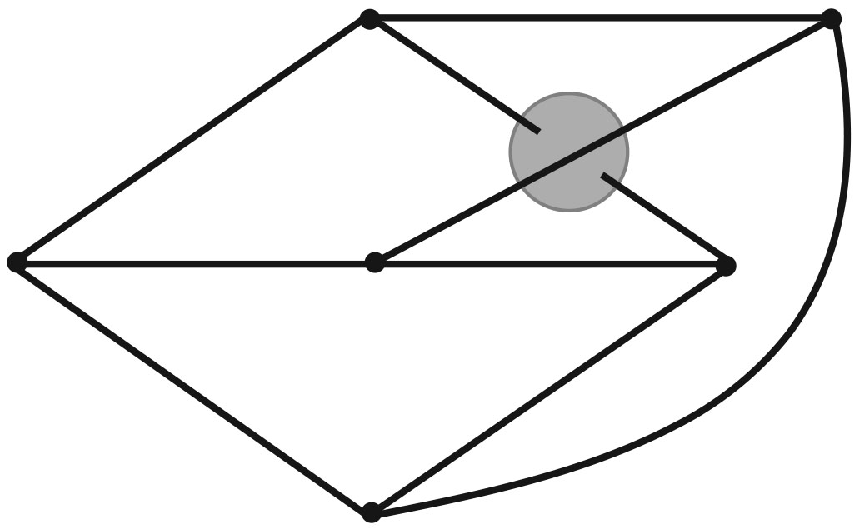}
\caption{}
\label{K33}
\end{figure}

In this section we develop methods to prove this conjecture for a certain family of non-planar cubic graphs. We will call a graph $G$ {\em near-planar} if it admits a planar projection with a single crossing, for example $K_{3,3}$ in figure \ref{K33}. We will view such graphs as $G=\Over\cup \overline G$ where $\overline G$ is a planar graph in a disk (the complement in $S^2$ of the shaded disk in figure \ref{K33}) with $4$ endpoints on the boundary circle. 

\begin{figure}[ht]
\includegraphics[height=2cm]{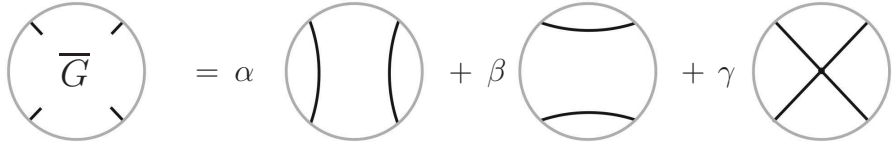}
{\large \put(-318,23){$\overline G$}}
\put(-267,25){$=$}
    \put(-252,25){$\alpha$}
    \put(-171,25){$+\; \, \beta$}
    \put(-83,25){$+ \:\, \gamma$}
\caption{A planar cubic graph in $C^{Q}_2$, $Q=(3-\sqrt{5})/2$, expressed as a linear combination of basis vectors.}
\label{element in C4}
\end{figure}

Considered as an element of the chromatic algebra $C^{Q}_2$, $\overline G$ may be expressed as a linear combination of the three basis elements without internal edges. The coefficients depend on $Q$; denote by ${\alpha}, {\beta}, {\gamma}$ their values at $Q=(3-\sqrt{5})/2$, figure \ref{element in C4}. In the following lemma the graph $\overline G$ outside a disk will be fixed throughout the proof, and (as in section \ref{notation section}) the  notation ${\circlSmall}_{\, Q}$ will stand for the evaluation of the flow polynomial at $Q$.

\begin{lem}[{\bf Golden identity for the flow polynomial of near-planar graphs}] \label{One crossing golden}  \sl 
Let $G=\Over\! \cup \overline{G}$ be a near-planar graph. Then the following identity holds for the evaluations of the flow polynomial of $G$ and of the two planar graphs obtained from $G$ by resolving the crossing: 
\begin{equation} \label{One crossing golden eq}
{\Over}_{\!\! z}\, =\, (-{\phi}^{-E}) \left( {\ZeroRes}_{\!\! w}^{\! 2} \, +\, {\phi}^{-1} {\ZeroRes}_{\!\! w}\,{\OneRes}_{\!\! w} \, +\, {\OneRes}_{\!\! w}^{\! 2} \right ),
\end{equation}
where $z=(5-\sqrt{5})/2,  w=(3-\sqrt{5})/2$.
\end{lem}

We remark that this lemma is motivated by D.W. Hall's version of the golden identity for constrained chromials \cite{Hall}.
The proof of lemma \ref{One crossing golden} relies on the following statement.

\begin{prop} \label{one crossing linear relation} \sl 
Fix $Q\neq (3\pm \sqrt{5})/2$, and let $G={\Over }\! \cup\overline G$ be a near-planar graph. Then 
$$
F^{}_G(Q)\, =\, {\Over}_{\!\! Q} \,
  =\, - \frac{1}{Q^2-3Q+1}\, {\ZeroRes}_{\!\! Q} \,  +\,  \frac{Q-1}{Q^2-3Q+1}\, {\Xfig}_{\!\! Q} \, - \frac{1}{Q^2-3Q+1} \, {\OneRes}_{\!\! Q}.$$
\end{prop}

The proof of proposition \ref{one crossing linear relation} amounts to checking that both sides have identical evaluations when paired up with the three basis elements $\ZeroRes$, $\OneRes$ and $\Xfig$ of the chromatic algebra ${\mathcal C}^Q_2$. (As discussed in section \ref{map}, the bilinear pairing $\langle .\, , \,  . \rangle_Q$ on $C_2^Q$ is non-degenerate  precisely for $Q\neq (3\pm \sqrt{5})/2$.)

{\em Proof of lemma \ref{One crossing golden}}.
Setting $Q$ to equal $z=(5-\sqrt{5})/2$ in proposition \ref{one crossing linear relation} and using the contraction-deletion rule, one has the following equality:
\begin{equation} \label{one crossing equation}
{\Over}_{\!\! z} \, = \, \\
 \frac{\phi}{2}\,  {\ZeroRes}_{\!\! z} \,  - \, \frac{{\phi}^{-1}}{2}\,  {\Ifig}_{\!\! z} \, +\,  \frac{1}{2} \,  {\OneRes}_{\!\! z}
\end{equation}

Equation  (\ref{One crossing golden eq}) is obtained by applying the golden identity (\ref{flow golden identity}) to the three planar cubic graphs on the right in (\ref{one crossing equation}).
\qed

It follows from lemma \ref{One crossing golden} that conjecture \ref{golden conjecture} for near-planar graphs is equivalent to the inequality
\begin{equation} \label{one crossing conj1}
{\Over}^{\! 2}_{\!\! w} \, \leq \, {\ZeroRes}^{\! 2}_{\!\! w} \, +\, {\OneRes}^{\! 2}_{\!\! w} \, + \, {\phi}^{-1} \, {\ZeroRes}_{\!\! w} {\OneRes}_{\!\! w} \!.
\end{equation}

\begin{remark} \rm
It is interesting to note that the  inequality (\ref{one crossing conj1}) can be restated in terms of the Yamada polynomial:
\begin{equation} \label{Yamada near planar}
 F^{}_{\Over}\!\! \left (\frac{3-\sqrt{5}}{2} \right ) \; \leq \; R^{}_{\Over}\!\! (e^{4{\pi} i/5})\; R^{}_{\Under}\!\! (e^{4{\pi} i/5}),
 \end{equation}

where the left-hand side is the evaluation of the flow polynomial  of the abstract near-planar graph, and the right-hand side is the product of the Yamada polynomials of two spatial ribbon graphs corresponding to the two possible crossings. The proof of the equivalence of (\ref{one crossing conj1}), (\ref{Yamada near planar}) consists of applying the Yamada polynomial skein relations to the crossings, and simplifying using the linear relation (\ref{linear Tutte1}).
\end{remark}
We are now in a position to give a reformulation of conjecture \ref{golden conjecture} for near-planar graphs.

\begin{lem} \label{one crossing conjecture} \sl
Let $\overline G$ be a cubic graph in a disk, with $4$ marked points on the boundary. Considered in the chromatic algebra $C^Q_2$, $Q=(3-\sqrt{5})/2$, $\overline G={\alpha} \ZeroRes+{\beta}\OneRes+{\gamma}\Xfig\!\!$. Conjecture \ref{golden conjecture} for near-planar graphs is equivalent to the  inequality
\begin{equation} \label{one crossing conjecture eq}
(1+3{\phi})\, {\alpha}\, {\beta} \; \leq \; {\gamma}\, ({\alpha}+{\beta}+{\gamma}).
\end{equation}
\end{lem}

{\em Proof.} 
The term being squared on the right-hand side of (\ref{conjugate inequality}) equals:
$$ F^{}_G(w)=\langle \Over,\overline G\rangle_{w}=\langle \Over, {\alpha} \ZeroRes+{\beta}\OneRes+{\gamma}\Xfig \rangle_{w}$$
$$=-{\phi}^{-1}({\alpha}+{\beta})+{\phi}^{-2}{\gamma}.$$
Analogous calculations for $\ZeroRes$ and $\OneRes$ in place of $\Over$ yield 
$$  -{\phi}^{-1}{\beta}+{\phi}^{-2}({\alpha}+{\gamma}), \; \;   -{\phi}^{-1}{\alpha}+{\phi}^{-2}({\beta}+{\gamma})
$$
respectively.
Equation (\ref{One crossing golden eq}) expresses the left-hand side in (\ref{conjugate inequality}) in terms of ${\ZeroRes}_{\!\! w}$ and ${\OneRes}_{\!\! w}$. Multiplying out the resulting expressions, (\ref{conjugate inequality}) is seen to be equivalent to (\ref{one crossing conjecture eq}). \qed

 Next we establish the  inequality  (\ref{one crossing conjecture eq})  for an infinite family of near-planar graphs.

\begin{figure}[ht]
\includegraphics[height=2cm]{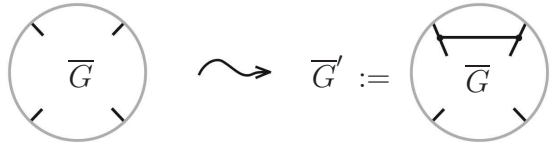}
 \put(-185,23){$\overline G$}
 \put(-93, 23){$\overline G'\, :=$}
    \put(-35,21){$\overline G$}
\caption{A modification of $\overline G$: addition of a {\em peripheral} edge (connecting two boundary edges of a cubic graph ${\overline G}$).}
\label{modification}
\end{figure}

\begin{figure}[ht]
\includegraphics[height=3cm]{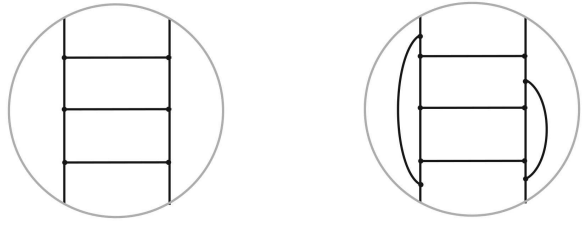}
\caption{Examples of graphs $\overline G$ in lemma \ref {family of graphs}. }
\label{examples figure}
\end{figure}

\begin{lem} \label{family of graphs} \sl
Conjecture \ref{golden conjecture} holds for the family of near-planar graphs $G={\Over} \! \cup \overline G$, where $\overline G$ is inductively built by addition of peripheral edges.
\end{lem}
Examples of graphs, considered in this lemma, are shown in figure \ref{examples figure}. The graph on the left, capped with ${\Over}\!\!$, is $K_{3,3}$. An analogous proof shows that the inequality (\ref{conjugate conjecture}) at the Galois conjugate values ${\phi}+1$, ${\phi}+2$ also holds for the same family of graphs.

\begin{figure}[ht]
\includegraphics[height=2cm]{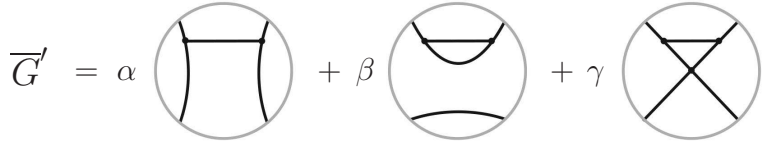}
{\large \put(-287,23){$\overline G'$}}
\put(-263,25){$=$}
    \put(-247,25){$\alpha$}
    \put(-171,25){$+\; \, \beta$}
    \put(-83,25){$+ \:\, \gamma$}
\caption{A peripheral edge added to $\overline G$.}
\label{modified element in C4}
\end{figure}

{\em Proof of lemma \ref{family of graphs}.}
It suffices to prove that the inequality (\ref{one crossing conjecture eq}) is preserved under addition of a peripheral edge. We will also show that whenever the coefficients ${\alpha}, {\beta}, {\gamma}$ are non-zero, they satisfy
\begin{equation} \label{signs}
{\rm sign}(\alpha)={\rm sign}(\beta)=-{\rm sign}(\gamma)=(-1)^{V/2},
\end{equation}

 where $V$ is the number of vertices of $\overline G$.
  The result of adding a peripheral edge to $\overline G={\alpha} \ZeroRes+{\beta}\OneRes+{\gamma}\Xfig\!\!$  is the graph $\overline G'$ shown in figure \ref{modified element in C4}. (There are $4$ possible peripheral edges; by symmetry between $\alpha$ and $\beta$ the proof below applies to each one.) 

Applications of the contraction-deletion rule give $$\overline G'=-{\alpha}\,  \ZeroRes +(-{\phi}{\beta} +{\gamma})\, \OneRes +({\alpha}-{\phi}^2{\gamma})\, \Xfig\!\!.$$

Denote the coefficients of $\overline G'$ by ${\alpha}', {\beta}', {\gamma}'$:
\begin{equation} \label{coefficient change}
{\alpha}' = -{\alpha},\; {\beta}'=-{\phi}{\beta} +{\gamma},\;  {\gamma}'={\alpha}-{\phi}^2{\gamma}.
\end{equation}

This gives an inductive proof of the statement (\ref{signs}). 
A direct calculation shows that the desired inequality for the coefficients  ${\alpha}', {\beta}', {\gamma}'$ of $\overline G'$,
\begin{equation} \label{inequality prime}
(1+3{\phi})\, {\alpha}'\, {\beta}' \; \leq \; {\gamma}'\, ({\alpha}'+{\beta}'+{\gamma}'),
\end{equation}

is equivalent to 
\begin{equation} \label{another inequality prime}
{\phi}^2\, {\alpha}\, {\beta} \; \leq \; {\gamma}\, ({\alpha}+{\beta}+{\gamma}).
\end{equation}

Since sign($\alpha$)$=$ sign($\beta)$, ${\phi}^2 \, {\alpha}\, {\beta}\leq(1+3{\phi})\, {\alpha}\, {\beta}$.Therefore the assumed inequality (\ref{one crossing conjecture eq}) for ${\alpha}, {\beta}, {\gamma}$ implies the inequality (\ref{inequality prime})  for ${\alpha}', {\beta}', {\gamma}'$.

\begin{figure}[ht]
\includegraphics[height=2.2cm]{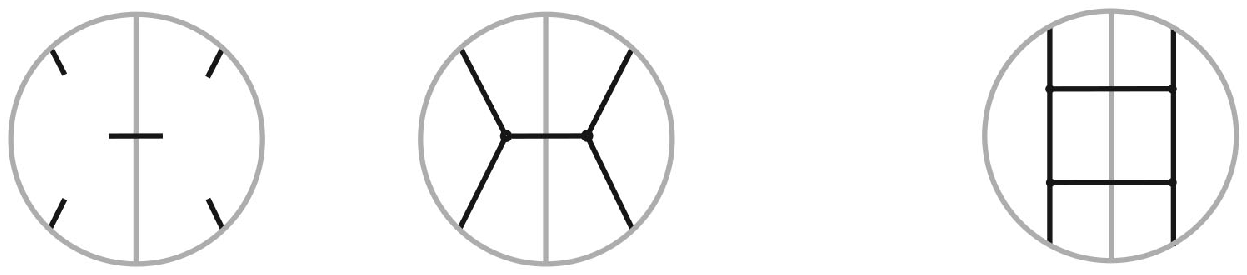}
\put(-302, 26){$\overline G=$}
{\small \put(-272,27){$H_1$}
\put(-238,27){$H_2$}}
    \put(-214,26){$={\lambda}$}
\caption{}
\label{tripod}
\end{figure}

Finally, consider the last statement in conjecture \ref{golden conjecture}. 
Note that the inequality ${\phi}^2 \, {\alpha}\, {\beta}\leq(1+3{\phi})\, {\alpha}\, {\beta}$ in the previous paragraph is strict precisely when both ${\alpha}, {\beta}$ are non-zero.
Suppose the inductive construction of the family of graphs in the statement of the lemma starts with $\Hfig$, and consider the first time a horizontal peripheral edge is added. Before this step, the graph $\overline G$ is of the form shown on the left in figure \ref{tripod}: there is a vertical line in the disk, intersecting the graph $\overline G$ in a single edge. The vector space of graphs in a disk with three boundary points, modulo the defining local relations of the chromatic algebra,  is $1$-dimensional. Hence it is clear that $\overline G = {\lambda} {\Hfig}$ in ${\mathcal C}_2$, and the coefficient $\beta$ (defined in figure) \ref{element in C4} is zero.
Then adding a horizontal edge gives a scalar multiple of the graph on the right in figure \ref{tripod}, resulting in non-zero coefficients ${\alpha}, {\beta}$. At this point, the graph $G={\Over}\!\!\cup \overline G$ is still planar, since one of the crossing strands may be drawn within the rectangle in the graph on the right in the figure. Since both coefficients ${\alpha}, {\beta}$ are non-zero, the addition of any new peripheral edge after this step makes the inequality (\ref{inequality prime}) strict. (And of course the graph becomes non-planar.) It follows from (\ref{coefficient change}) that all further additions of edges increase ${\alpha}, {\beta}$ in absolute value, so the difference of the two sides of the inequality  (\ref{inequality prime}) continues strictly increasing. Thus the inequality detects planarity in this family of graphs.
\qed

\section{Exponential growth of the number of chromatic polynomials}  \label{number of chromatic polynomials}

In this section we prove theorem \ref{ChromaticPolynomials thm}, working dually with the flow polynomial of planar cubic graphs. Consider ${\mathcal C}^Q_{1,3}$, the vector space consisting of  ${\mathbb C}$-linear combinations of planar cubic graphs in a rectangle with one marked point at the bottom and three marked points at the top of the rectangle, modulo the relations in figure \ref{chromatic relations}. (This notion was also discussed at the end of  section \ref{flow category subsection}.)
The loop value is $Q-1$, and for the remainder of this section we fix $Q=(3-\sqrt{5})/2$. As discussed in section \ref{map}, at this value of $Q$ there is an additional local relation 
\begin{equation} \label{linear Tutte000}
\Ifig\, =\, -{\phi}\, \ZeroRes - {\phi}^{2}\,  \OneRes.
\end{equation}
${\mathcal C}^{Q}_{1,3}$ is a module over the chromatic algebra; the action is by vertical concatenation, matching three marked points on the boundary.

\begin{figure}[h]
\centering
\includegraphics[height=2.2cm]{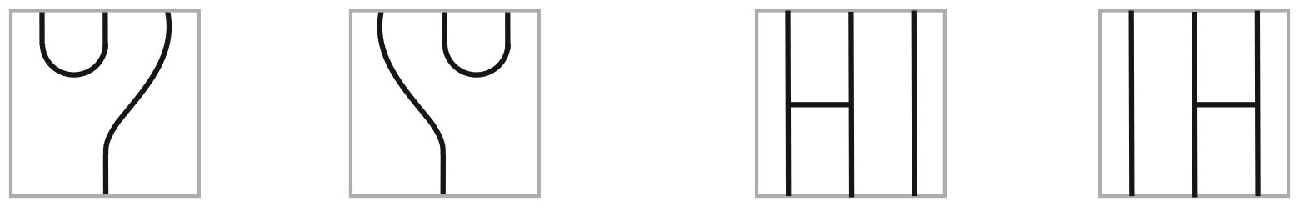}
{\small     \put(-356,27){$e_1 =$}
    \put(-270,27){$e_2 =$}
     \put(-170,27){$A = $}
    \put(-82,27){$B = $}
    }  
    \caption{Elements  $e_1, e_2\in {\mathcal C}^{Q}_{1,3}$ and $A,B\in {\mathcal C}^Q_3$}
\label{TL0 figure}
\end{figure}

\begin{figure}[h]
\centering
\includegraphics[height=2.2cm]{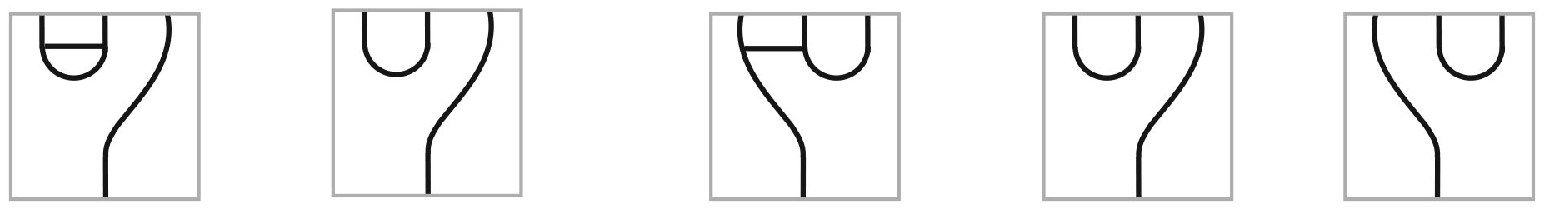}
{\small     \put(-416,27){$A e_1\! =$}
    \put(-336,27){$ =-{\phi}$}
    \put(-253,27){$,$}
     \put(-243,27){$Ae_2 = $}
    \put(-159,27){$=-{\phi} $}
       \put(-80,27){$-{\phi}^2 $}
    }  
    \caption{The action of $A$}
\label{TL1 figure}
\end{figure}

\begin{figure}[h]
\centering
\includegraphics[height=2.2cm]{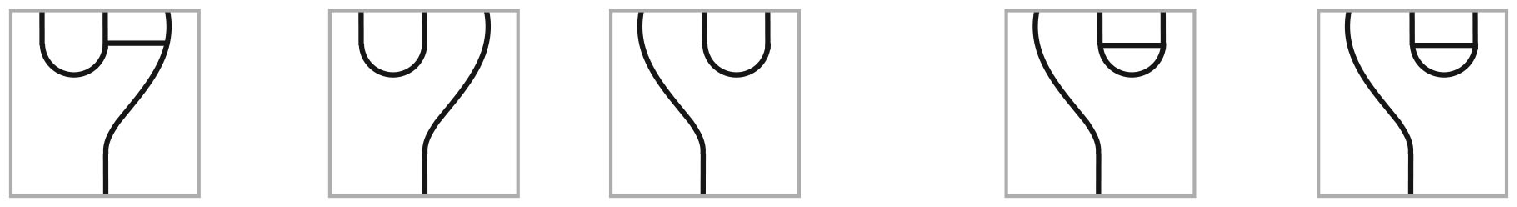}
{\small     \put(-416,27){$B e_1\! =$}
    \put(-337,27){$ =\! -{\phi}^2$}
    \put(-253,27){$-{\phi}$}
     \put(-181,27){$,$}
    \put(-165,27){$B e_2 =$}
       \put(-83,27){$=-{\phi} $}
    }  
    \caption{The action of $B$}
\label{TL2 figure}
\end{figure}

For a generic $Q$, the vector space ${\mathcal C}^Q_{1,3}$ is $3$-dimensional. At the specified value $Q=(3-\sqrt{5})/2$, the additional relation (\ref{linear Tutte000}) reduces the dimension to two. Consider its basis $e_1, e_2$, shown on the left in figure \ref{TL0 figure}, and
let $A, B$ be elements of  ${\mathcal C}_3$ shown on the right in the same figure. The action of $A, B$ on $e_1, e_2$ is calculated in figures \ref{TL1 figure}, \ref{TL2 figure}.

In other words, with respect to the chosen basis, $A$ and $B$ are represented by the matrices $$A=-{\phi}\begin{bmatrix} 1 & 1 \\0 & {\phi}
\end{bmatrix}, B=-{\phi}\begin{bmatrix} {\phi} & 0 \\1 & 1
\end{bmatrix}.$$

The squares of these matrices are given by
$$
A^2={\phi}^2\begin{bmatrix} 1 & {\phi}^2 \\0 & {\phi}^2
\end{bmatrix}, B^2={\phi}^2\begin{bmatrix} {\phi}^2 & 0 \\{\phi}^2 & 1
\end{bmatrix}.$$

It follows from the ping pong lemma for semigroups that $A^2, B^2$ generate a free semigroup in the chromatic algebra ${\mathcal C}^Q_3$. (Consider vectors $v$ with positive components with respect to the basis $e_1, e_2$, and let 
$u=A^2 v, w=B^2 v$. Then the components of $u, w$ satisfy $u_1>u_2, w_1<w_2$.)
Therefore words $\{ W\}$  of length $n$ in $A^2, B^2$ represent $2^n$ distinct elements in ${\mathcal C}^Q_3$. 

\begin{figure}[h]
\centering
\includegraphics[height=2.4cm]{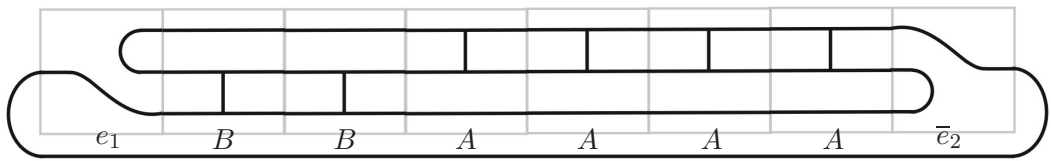}
{\small     \put(-368,12){$e_1$}
    \put(-278,10){$B$}
      \put(-324,10){$B$}
     \put(-186,10){$A $}
        \put(-232,10){$A $}
     \put(-139,10){$A $}
    \put(-93,10){$A $}
        \put(-50,12){$\overline e_2 $}
    }  
    \caption{The flow polynomial of the pictured graph at $(3-\sqrt{5})/2$ equals $\langle  A^4B^2e_1, e_2 \rangle$. (The multiplication is represented horizontally, rather than vertically, to conserve space.)}
\label{TL3 figure}
\end{figure}

Given two elements $a, b\in {\mathcal C}^Q_{1,3}$, the product $a \cdot \overline b$ (where $\overline b \in {\mathcal C}^Q_{3,1}$ is obtained from $b$ by reflection in a horizontal line) is an element of ${\mathcal C}^Q_1$. Here the product is given by vertical concatenation, matching three marked boundary points. ${\mathcal C}_1$ is  $1$-dimensional, and the trace ${\mathcal C}^Q\longrightarrow {\mathbb C}$ is an isomorphism.  Denote $\langle a, b \rangle: = {\rm tr}(a \cdot \overline b)$. 
Note that the Gram matrix $\langle e_i, e_j \rangle_{i,j\in \{1, 2\} }$ is non-degenerate at $Q=(3-\sqrt{5})/2$. Therefore if two words $W_1, W_2$ in $A^2, B^2$  are not equal as $2 \times 2$ matrices, then $\langle W_1 e_i, e_j\rangle \neq \langle W_2 e_i, e_j \rangle$ for some $i,j\in \{1, 2\}$.
Words in $A^2, B^2$ are represented geometrically by planar cubic graphs $G$ in a rectangle, and (as discussed in section \ref{chromatic algebra subsection}) $\langle W e_i, e_j\rangle$ is the flow polynomial of the graph obtained by gluing $e_i$ on the bottom of $G$, $\overline e_j$ on top, and then taking the trace, i.e. connecting the endpoints by an arc in the plane, figure \ref{TL3 figure}.

Since words of length $n$ in $A^2, B^2$ give $2^n$ distinct matrices, there are at least $2^{n/4}$ distinct values for one of the matrix entries. This translates into at least $2^{n/16}$ values $F^{}_G((3-\sqrt{5})/2)$ of the flow polynomial of planar trivalent graphs $G$ with $n$ vertices.
\qed

{\bf Remark:} After the paper was written, we discovered that one may also prove exponential growth at $Q=4$, observing that $A^2e_1=4e_1,A^2e_2=e_1+e_2, B^2e_1=e_1+e_2, B^2e_2=4e_2$, and completing the argument as above.

\section{Appendix: A proof of the golden identity for planar cubic graphs} \label{appendix}

This section gives a variation of the proof \cite{FK} of the Tutte golden identity \cite{T2} for the flow polynomial of planar cubic graphs:
\begin{equation} \label{flow golden identity again}
F^{}_G((5-\sqrt{5})/2) =  (-{\phi})^{-E}\, F^{}_G((3-\sqrt{5})/{2})^2,
\end{equation}
We include a proof since this identity underlies several results in this paper. The version of the proof presented here has more of a computational flavor; it might be useful in numerical investigation of whether there are identities at other parameter values.

Denote $z:=(5-\sqrt{5})/2,  w:=(3-\sqrt{5})/2$. Let $\overline G$ be a graph in the disk with $4$ marked points on the boundary, and consider its coefficients with respect to the usual basis of $C^Q_2$ (as in figure \ref{element in C4}), 
$$\overline G \, =\, {\alpha}_Q\; {\ZeroRes}_{\!\! Q} +{\beta}_Q\; {\OneRes}_{\!\! Q} + {\gamma}_Q\; {\Xfig}_{\!\! Q}$$
for $Q=z, w$. Flips  $\Ifig \leftrightarrow \Hfig$ acts transitively on connected planar cubic graphs with a fixed number of  vertices, and it suffices to prove that the golden identity for $\ZeroRes\!\!$, $\OneRes\!\!$, and $\Ifig\!\!$ imply the golden identity for $\Hfig$.

Consider the evaluations:
 $$\langle \ZeroRes, \overline G\rangle_Q \; =\;  {\alpha}_Q (Q-1)^2+{\beta}_Q (Q-1)+{\gamma}_Q (Q-1)^2 ,$$ 
$$ \langle \OneRes, \overline G\rangle_Q \;  =\;  {\alpha}_Q (Q-1)+{\beta}_Q (Q-1)^2+{\gamma}_Q (Q-1)^2 ,$$ 
$$ \langle \Ifig, \overline G\rangle_Q \;  =\;  {\alpha}_Q (Q-1)(Q-2)+{\gamma}_Q (Q-1)(Q-2)^2,$$ 
$$ \langle \Hfig, \overline G\rangle_Q \; =\; {\beta}_Q (Q-1)(Q-2)+{\gamma}_Q (Q-1)(Q-2)^2 .$$

The golden identity for $G= \overline G$, capped off  with $\ZeroRes, \OneRes, \Ifig, \Hfig$, respectively, is equivalent to quadratic equations on the coefficients at $z$ and $  w$:
\begin{equation} \tag{$A$} {\alpha}_z {\phi}^{-4}+{\beta}_z {\phi}^{-2}+{\gamma}_z {\phi}^{-4}\; -\; \left ({\alpha}_w {\phi}^{-2}-{\beta}_w {\phi}^{-1}+{\gamma}_w {\phi}^{-2} \right )^2 \; = \; 0,
\end{equation}
\begin{equation} \tag{$B$}{\alpha}_z {\phi}^{-2}+{\beta}_z {\phi}^{-4}+{\gamma}_z {\phi}^{-4}\; -\; \left (-{\alpha}_w {\phi}^{-1}+{\beta}_w {\phi}^{-2}+{\gamma}_w {\phi}^{-2}\right )^2\; = \; 0,
\end{equation}
\begin{equation} \tag{$C$}-{\alpha}_z {\phi}^{-3}+{\gamma}_z {\phi}^{-4}\; -\; \left ({\alpha}_w -{\gamma}_w {\phi}\right )^2\; = \; 0,
\end{equation}
\begin{equation} \tag{$D$}-{\beta}_z {\phi}^{-3}+{\gamma}_z {\phi}^{-4}\; -\; \left ({\beta}_w -{\gamma}_w {\phi}\right )^2\; = \; 0.
\end{equation}

The left-hand sides of these four equations satisfy the relation $({\phi}^{-2}-{\phi}^{-4}) (A-B) =C-D$, so the validity of any three of them implies the fourth.
The proof of (\ref{flow golden identity again}) is completed by comparing the loop values: ${\circlblSmall}_{\;\, z}={\circlblSmall}_{\;\, w}^{\;\, 2}$. 
\qed

{\bf Acknowledgements.} We would like to thank Gordon Royle for many discussions, and also for sharing with us  numerical data on graph polynomials. We also thank Kyle Miller for helpful comments. 

V. Krushkal was supported in part by NSF grant DMS-1612159. Ian Agol was supported by NSF grant DMS-1406301, a Simons Investigator Award, and the Institute for Advanced Study.


\begin{thebibliography}{10}

\bibitem{AK} I. Agol and V. Krushkal,
{\em Tutte relations, TQFT, and planarity of cubic graphs}, 
Illinois J. Math. 60 (2016), no. 1, 273-288. 


\bibitem{Beraha} S. Beraha, {\em Infinite non-trivial families of maps and chromials}, Thesis, Johns Hopkins University, 1975.

\bibitem{BL} G.D. Birkhoff and D.C. Lewis, {\em  Chromatic polynomials}, Trans. Amer. Math. Soc. 60 (1946), 355-451.

\bibitem{CS} B. Chen and R. Stanley, {\em Orientations, lattice polytopes, and group arrangements II: Modular and integral flow polynomials of graphs}, Graphs Combin. 28 (2012), no. 6, 751-779. 


\bibitem{Diestel}  R. Diestel, Graph theory. Graduate Texts in Mathematics, 173. Springer-Verlag, Berlin, 2005.

\bibitem{FK} P. Fendley and V.  Krushkal, {\em Tutte chromatic identities from the Temperley-Lieb algebra}, Geom. Topol. 13 (2009), 709-741.

\bibitem{FK2} P. Fendley and V.  Krushkal, {\em Link invariants, the chromatic polynomial and the Potts model}, Adv. Theor. Math. Phys.  14 (2010) 2, 507--540.

\bibitem{Hall} D.W. Hall, {\em On golden identities for constrained chromials},  J. Combinatorial Theory Ser. B 11 (1971), 287-298.

\bibitem{Jaeger} F. Jaeger,  {\em On the Penrose number of cubic diagrams},
Graph colouring and variations. 
Discrete Math. 74 (1989), 85-97.

\bibitem{Jo} V.F.R.~Jones,
Subfactors and knots. 
CBMS Regional Conference Series in Mathematics, 80. Published for the Conference Board of the Mathematical Sciences, Washington, DC; American Mathematical Society, Providence, RI, 1991.

\bibitem{KV} L.H. Kauffman and P. Vogel, {\em Link polynomials and a graphical calculus}, J. Knot Theory Ramifications 1 (1992),
59-104.

\bibitem{MW} P. Martin and D. Woodcock, 
{\em The partition algebras and a new deformation of the Schur algebras.} 
J. Algebra 203 (1998), no. 1, 91-124. 

\bibitem{MPS} S. Morrison, E. Peters and N. Snyder, {\em  Knot polynomial identities and quantum group coincidences}, Quantum Topol. 2 (2011), 101-156.

\bibitem{Penrose} R. Penrose, {\em Applications of negative dimensional tensors}, 1971 Combinatorial Mathematics and its Applications (Proc. Conf., Oxford, 1969) pp. 221-244. Academic Press, London.

\bibitem{RT} N.Yu. Reshetikhin and V.G. Turaev, {\em  Ribbon graphs and their invariants derived from quantum groups}, Comm. Math. Phys. 127 (1990), 1-26.

\bibitem{Treumann} D. Treumann, {\em How many chromatic polynomials of planar maps are there?}, https://mathoverflow.net/q/241734

\bibitem{TZ} D. Treumann and E. Zaslow, {\em Cubic Planar Graphs and Legendrian Surface Theory}, 
arXiv:1609.04892


\bibitem{T63} W.T. Tutte, {\em A census of planar maps}, 
Canad. J. Math. 15 1963 249-271. 



\bibitem{T1} W.T. Tutte, {\em On chromatic polynomials and the golden ratio}, J. Combinatorial Theory 9 (1970), 289-296.


\bibitem{T2} W.T.~Tutte, {\em More about chromatic polynomials and the golden ratio}, Combinatorial Structures and their Applications, 439-453 (Proc. Calgary Internat. Conf., Calgary, Alta., 1969)



\bibitem{T3} W.T. Tutte, {\em A contribution to the theory of chromatic polynomials},
Canadian J. Math. 6, (1954), 80-91.

\bibitem{Walker} K. Walker, Deletion-contraction relations, hard hexagons, and the shadow world, Talk at IPAM, 2007.
Available at http://canyon23.net/math/talks/

\bibitem{We} H.Wenzl,
{\em On a sequence of projections},  C. R. Math. Rep. Acad. Sci. Canada  9  (1987), 5-9.

\bibitem{Yamada} S. Yamada, {\em An invariant of spatial graphs}, J. Graph Theory 13 (1989), 537-
551.




\end{thebibliography}
\end{document}